\documentstyle[12pt]{article}
\def\txt#1{{\textstyle{#1}}}
\def\hf{{\textstyle{1\over2}}}
\def\a{\alpha}
\def\d{{\,\rm d}}
\def\DJ{\leavevmode\setbox0=\hbox{D}\kern0pt\rlap
 {\kern.04em\raise.188\ht0\hbox{-}}D}
\def\e{\varepsilon}
\def\f{\varphi}
\def\G{\Gamma}

\def\k{\kappa}
\def\s{\sigma}

\def\={\;=\;}
\def\zx{\zeta(\hf+ix)}
\def\zt{\zeta(\hf+it)}

\def\D{\Delta}
  
\def\R{\Re{\rm e}\,}  \def\s{\sigma}
\def\z{\zeta}

\def\no{\noindent} 
\def\H{H_j^3({\txt{1\over2}})}  \def\={\,=\,}
\def\hf{{\textstyle{1\over2}}}
\def\txt#1{{\textstyle{#1}}}
\def\f{\varphi}
\def\Z{{\cal Z}}

\font\tenmsb=msbm10
\font\sevenmsb=msbm7
\font\fivemsb=msbm5
\newfam\msbfam
\textfont\msbfam=\tenmsb
\scriptfont\msbfam=\sevenmsb
\scriptscriptfont\msbfam=\fivemsb
\def\Bbb#1{{\fam\msbfam #1}}

\def \NN {\Bbb N}
\def \CC {\Bbb C}

\def \ZZ {\Bbb Z}

\title{\sc The Laplace and Mellin transforms of
powers of the Riemann zeta-function}
\medskip
\author{Aleksandar Ivi\'c}
\begin{document}
\maketitle
\begin{abstract}

\bigskip
This paper gives a survey of known results concerning
the Laplace transform
$$
L_k(s) \;:=\; \int_0^\infty |\zx|^{2k}{\rm e}^{-sx}\d x
\qquad(k \in \NN,\, \R s > 0),
$$
and the (modified) Mellin transform
$$
\Z_k(s) \;:= \;\int_1^\infty|\zx|^{2k}x^{-s}\d x\qquad(k\in\NN),
$$
where the integral is absolutely
convergent for $\R s \ge c(k) > 1$.
Also some new results on these  integral transforms
of $|\zx|^{2k}$ are given, which have important
connections with power moments of the Riemann zeta-function $\z(s)$.
\end{abstract}
\section{ Introduction}

A central place in Analytic Number Theory is occupied by the Riemann
zeta-function $\zeta(s)$, defined  by the representations
$$
\zeta(s) = \sum_{n=1}^\infty n^{-s} =
\prod_{p\,:\,{\rm{prime}}} {(1 - p^{-s})}^{-1}
\qquad(s = \s+it,\;\s>1), \eqno{(1.1)}
$$
and otherwise by analytic continuation. It admits meromorphic
continuation to the whole complex plane, its only singularity
being the simple pole $s = 1$ with residue 1. For general
information on $\zeta(s)$ the reader is referred to the
monographs [5] and [30]. From the functional equation
$$
\zeta(s) = \chi(s)\zeta(1 - s), \quad\chi(s) = 2^s\pi^{s-1}
\sin\bigl({{\pi s}\over 2}\bigr)\Gamma(1 - s), \eqno{(1.2)}
$$
which is valid for any complex $s$, it follows that $\zeta(s)$ has
zeros at $s = -2,-4,\ldots$~. These zeros are traditionally called
the ``trivial'' zeros of $\zeta(s)$, to distinguish them from the
complex zeros of $\zeta(s)$, of which the smallest ones (in
absolute value) are ${1 \over 2} \pm 14.134725\ldots i$. It is
well-known that all complex zeros of  $\zeta(s)$  lie in the
so-called ``critical strip'' $0 < \sigma = \R s < 1$, and if
$N(T)$ denotes the number of zeros $\rho = \beta + i\gamma$
($\beta, \gamma$ real) of  $\zeta(s)$ for which $0 < \gamma \le
T$, then
$$
N(T) = {T\over {2\pi}}\log\bigl({T\over {2\pi}}\bigr) - {T\over {2\pi}}
+ {7\over 8} + S(T) + O({1\over T}) \eqno{(1.3)}
$$
with
$$
S(T) = {1\over \pi}\arg \zeta(\txt{1\over 2} + iT) = O(\log T). \eqno{(1.4)}
$$
Here $S(T)$ is obtained by continuous variation along the straight lines
joining $2, 2 + iT$, ${1\over2}+iT$, starting with the value $0$; if $T$ is
the ordinate of a zero, let $S(T)=S(T+0)$. This is the
so-called Riemann--von Mangoldt formula. {\it The Riemann hypothesis} (RH)
 is the
conjecture, stated by B. Riemann in his epoch-making memoir [29], that
``{\it very likely all complex zeros of  $\zeta(s)$
have real parts equal to $1\over2$}". For this reason the line
$\sigma = {1\over2}$ is called the ``critical line'' in the theory of
$\zeta(s)$. The RH is undoubtedly one of the most celebrated and difficult
open problems in whole Mathematics. Its proof (or disproof) would
have very important consequences in multiplicative number theory,
especially in problems involving the distribution of primes. It would also
very likely lead to generalizations to many other zeta-functions
(Dirichlet series) sharing similar properties with  $\zeta(s)$.

\bigskip
The aim of this paper is to present known results on the Laplace
transform
$$
L_k(s) \;:=\; \int_0^\infty |\zx|^{2k}{\rm e}^{-sx}\d x
\qquad(k \in \NN,\,\s = \R s > 0),\eqno(1.5)
$$
and the (modified) Mellin transform
$$
\Z_k(s) \;:=\; \int_1^\infty|\zx|^{2k}x^{-s}\d x\qquad(k\in\NN),\eqno(1.6)
$$
where the integral in (1.6) is absolutely convergent for $\R s \ge
c(k) > 1$, and also to present some new results. The term ``modified"
refers to the fact that the lower bound of integration in (1.6) is 1, and not
0 which is usual for Mellin transforms, and also we have $x^{-s}$
instead of $x^{s-1}$. These modifications are technical, as the lower
bound 1 dispenses with convergence problems which may occur when $|s|$
is large.

\medskip

Apart from the
distribution of complex zeros of $\z(s)$, there are several
central topics in zeta-function theory. One of them is doubtlessly
the evaluation (or estimation) of power moments of $|\z(\s+it)|$, that is,
integrals of the form $\int_0^T|\z(\s+it)|^{2k}\d t\;(0 < \s <
1)$, where $\s$ is assumed to be fixed.
The most important case, in view of the functional equation
(1.2), is when $\s = 1/2$ and also when $k\in\NN$, although the
case when $k$ is not an integer is also of interest (see e.g., [6, Chapter 6],
and in general for power moments of the zeta-function see [5], [6] and
[30]). Thus the most important object of study involving power
moments of $|\zt|$ is
$$
I_k(T) \;:=\; \int_0^T|\zt|^{2k}\d t\qquad(k \in \NN). \eqno(1.7)
$$
One trivially has
$$
I_k(T) \le {\rm e}\int_0^\infty|\zt|^{2k}
{\rm e}^{-t/T}\d t = {\rm e}L_k\left({1\over T}\right).
\eqno(1.8)
$$
Therefore any nontrivial bound of the form
$$
L_k(\s) \;\ll_\e\; \left({1\over\s}\right)^{c_k+\e}\qquad(\s\to 0+,\,
c_k \ge 1)\eqno(1.9)
$$
gives, in view of (1.8) ($\s = 1/T$), the bound
$$
I_k(T) \;\ll_\e\; T^{c_k+\e}.\eqno(1.10)
$$
Conversely, if (1.10) holds, then we obtain (1.9) from the identity
$$
L_k\left({1\over T}\right) \= {1\over T}\,\int_0^\infty I_k(t)
{\rm e}^{-t/T}\d t,
$$
which is easily established by integration by parts.
We note that here and later $\e$ denotes arbitrarily small constants,
not necessarily the same ones at each occurrence. The symbol $f(x) \ll g(x)$ (same
as $f(x)=O(g(x))$) means that $|f(x)| \le Cg(x)$ for some $C>0$ and $x\ge x_0$, while
$f(x) \ll_{a,b,\ldots} g(x)$ means that the $\ll$--constant
 depends on $a,b,\,\ldots\;$.

\medskip
One of the possible applications of ${\cal Z}_k(s)$ consists of the
following. If
$$F(s) \= \int_0^\infty f(x)x^{s-1}\d x
$$
is the (classical) Mellin transform of $f(x)$,
then by (4.2) obtains, for suitable $c > 1$,
$$
\int_1^\infty f\left({x\over T}\right)|\zx|^{2k}\d x
$$
$$=
\int_1^\infty{1\over2\pi i} \int\limits_{(c)}F(s)\left({T\over
x}\right)^s\d s|\zx|^{2k}\d x = {1\over2\pi
i}\int_{(c)}F(s)T^s{\cal Z}_k(s)\d s,\eqno(1.11)
$$
where as usual $\int_{(c)}$ denotes integration over the line $\R s = c$.
If $f(x) \in C^\infty$ is a nonnegative function of compact support
such that $f(x) = 1$ for $1 \le x \le 2$, then $F(s)$ is entire of
fast decay, and (1.11) (with $c = 1+\e$) yields a weak form of the
$2k$--th moment for $|\zt|$, namely
$$
\int_0^T\vert\zeta({\txt{1\over 2}} + it)\vert^{2k}\d  t\;
\ll_{k,\e} \; T^{1+\e},\eqno(1.12)
$$
provided that ${\cal Z}_k(s)$ has analytic continuation to the
half-plane $\s > 1$, where it is regular and of polynomial growth
in $|t|$. Conversely, if (1.12) holds, then integrating by parts
it is seen that ${\cal Z}_k(s)$  is regular for $\s > 1$
(see Theorem 2 below for the precise statement).

\medskip

Note that, by a change of variables $x = {\rm e}^t,\,z = s-1$, (1.5) becomes
$$
\int_0^\infty |\z(\hf + i{\rm e}^t)|^{2k}{\rm e}^{-zt}\d t\qquad
(\R z > 0),
$$
which is the Laplace transform of $|\z(\hf + i{\rm e}^t)|^{2k}$. Indeed, it
is well-known that the Laplace and Mellin transforms are closely connected, as
(by a change of variable) both of them can be regarded as special cases of
Fourier transforms, and their theory built from the theory of Fourier transforms.

\bigskip
\section{Laplace transforms of powers of the zeta-function}

Laplace transforms play an important r\^ole in analytic number
theory. Of special interest in the theory of the Riemann
zeta-function $\z(s)$ are the Laplace transforms (1.5), namely the functions
$$
L_k(s) \;:=\; \int_0^\infty |\zx|^{2k}{\rm }{\rm e}^{-sx}\d x\qquad(k \in
\NN,\,\R s > 0).
$$
E.C. Titchmarsh's well-known monograph [30, Chapter 7] gives a
discussion of $L_k(s)$ when $s = \s$ is real and $\s \to 0+$,
especially detailed in the cases $k=1$ and $k=2$. Indeed, a
classical result of H. Kober [24] says that, as $\s \to 0+$,
$$
L_1(2\s) = {\gamma-\log(4\pi\s)\over2\sin\s} +
 \sum_{n=0}^Nc_n\s^n + O_N(\s^{N+1})
\eqno(2.1)
$$
for any given integer $N \ge 1$, where the $c_n$'s are effectively
computable constants and $\gamma = -\G'(1) = 0.577\ldots\,$ is
Euler's constant. For complex values of $s$ the function $L_1(s)$
was studied by F.V. Atkinson [1], and more recently by M. Jutila
[22], who noted that Atkinson's argument gives
$$
L_1(s) = -i{\rm e}^{{1\over2}is}(\log(2\pi)-\gamma + ({\pi\over2}-s)i) +
2\pi {\rm e}^{-{1\over2}is}\sum_{n=1}^\infty d(n)\exp(-2\pi in{\rm e}^{-is}) +
\lambda_1(s) \eqno(2.2)
$$
in the strip $0 < \R s < \pi$, where the function $\lambda_1(s)$
is holomorphic in the strip $|\R s| < \pi$. Moreover, in any strip
$|\R s| \le \theta$ with $0 < \theta < \pi$, we have
$$
\lambda_1(s) \;\ll_\theta \;(|s|+1)^{-1}.
$$
In [21] M. Jutila gave a discussion on the application of Laplace
transforms to the evaluation of sums of coefficients of certain
Dirichlet series.

For $L_2(\s)$ F.V. Atkinson [2] obtained the asymptotic formula
$$
L_2(\s) = {1\over\s}\left(A\log^4{1\over\s} + B\log^3{1\over\s} +
C\log^2{1\over\s} + D\log {1\over\s} + E\right) +
\lambda_2(\s),\eqno(2.3)
$$
where $\s \to 0+$,
$$
A = {1\over2\pi^2},\,B =\pi^{-2}(2\log(2\pi) - 6\gamma +
24\z'(2)\pi^{-2})
$$
and
$$
 \lambda_2(\s) \;\ll_\e\;\left({1\over\s}\right)^{{13\over14}+\e}.
 \eqno(2.4)
$$
He also indicated how, by the use of estimates for Kloosterman
sums, one can improve the exponent ${13\over14}$ in (2.4) to
${8\over9}$. This is of historical interest, since it is one of
the first instances of an application of Kloosterman sums to analytic
number theory. Atkinson in fact showed that ($\s = \R s > 0$)
$$
L_2(s) \;=\;4\pi {\rm e}^{-{1\over2}s}\sum_{n=1}^\infty d_4(n) K_0(4\pi
i\sqrt{n}{\rm e}^{-{1\over2}s}) + \phi(s),\eqno(2.5)
$$
where $d_4(n)$ is the divisor function generated by $\z^4(s)$,
$K_0$ is the Bessel function, and the series in (2.5) as well as
the function $\phi(s)$ are both analytic in the region $|s| <
\pi$. When $s = \s \to 0+$ one can use the asymptotic formula
$$
K_0(z) \= \hf\sqrt{\pi}z^{-1/2}{\rm e}^{-z}\left(1 - 8z^{-1} +
O(|z|^{-2})\right) \quad\left(|\arg z| < \theta < {3\pi\over2},\,|z|
\ge 1\right)
$$
and then, by delicate analysis, one can deduce (2.3)--(2.4) from
(2.5).

\medskip
The author [7] gave explicit, albeit complicated expressions for
the remaining coefficients $C,D$ and $E$ in (2.3). More
importantly, he applied a result on the fourth moment of $|\zt|$,
obtained jointly with Y. Motohashi [19], [20] (see also [28]), to
establish that
$$
\lambda_2(\s) \;\ll\; \s^{-1/2}\qquad(\s\to 0+),\eqno(2.6)
$$
which in view of Theorem 1 below is best possible.

\medskip
For $k \ge 3$ not much is known about $L_k(s)$, even when $s = \s
\to 0+$. This is not surprising, since not much is known about
upper bounds for the integral $I_k(T)$ (see (1.7)).
For a discussion on $I_k(T)$ the reader is referred to the
author's monographs [5] and [6].

\medskip

One can consider $L_2(s)$, where $s$ is a
complex variable, and  prove a result analogous to (2.2), valid
in a certain region in $\CC$. This was achieved by the
author in [10]. The main tools that were used are the
results and methods from spectral theory, by which recently much advance has
been made in connection with $I_2(T)$ (cf. (1.7); see [6]--[20], and [26]--[28]
for some of the relevant papers on this subject). For a competent and
extensive account of spectral theory the reader is referred to Y.
Motohashi's monograph [28].

\smallskip

We shall state here very briefly  the necessary notation involving the
spectral theory of the non-Euclidean Laplacian. As usual
$\,\{\lambda_j = \kappa_j^2 + {1\over4}\} \,\cup\, \{0\}\,$ will
denote the discrete spectrum of the non-Euclidean Laplacian acting
on $\,SL(2,\ZZ)\,$--automorphic forms, and $\a_j =
|\rho_j(1)|^2(\cosh\pi\kappa_j)^{-1}$, where $\rho_j(1)$ is the
first Fourier coefficient of the Maass wave form corresponding to
the eigenvalue $\lambda_j$ to which the Hecke series $H_j(s)$ is
attached. We note that
$$
\sum_{\k_j\le K}\a_jH_j^3(\hf) \ll K^2\log^CK \qquad(C >
0).\eqno(2.7)
$$
Our result on $L_2(s)$, proved in [10], is the following

\medskip
THEOREM 1. {\it Let $0 \le \phi < {\pi\over2}$ be given. Then for
$0 < |s| \le 1$ and $|\arg s| \le \phi$ we have}
$$
L_2(s) = {1\over s}\left(A\log^4{1\over s} + B\log^3{1\over s} +
C\log^2{1\over s} + D\log{1\over s} + E\right)
$$
$$
+\, s^{-{1\over2}}\left\{\sum_{j=1}^\infty \a_j\H\Bigl(
s^{-i\k_j}R(\k_j)\G(\hf + i\k_j) + s^{i\k_j}R(-\k_j)\G(\hf -
i\k_j) \Bigr)\right\}  + G_2(s),
\eqno(2.8)
$$
{\it where}
$$
R(y) \;:=\; \sqrt{{\pi\over2}}{\Bigl(2^{-iy}{\G({1\over4} -
{i\over2}y) \over\G({1\over4} +
{i\over2}y)}\Bigr)}^3\G(2iy)\cosh(\pi y) \eqno(2.9)
$$
{\it and in the above region $G_2(s)$ is a regular function
satisfying }($C > 0$ {\it is a suitable constant})
$$
G_2(s) \ll |s|^{-1/2}\exp\left\{
-{C\log(|s|^{-1}+20)\over(\log\log(|s|^{-1}+20))^{2/3}
(\log\log\log(|s|^{-1}+20))^{1/3}}\right\}. \eqno(2.10)
$$

{\bf Remark 1}. The constants $A,\,B,\,C,\,D,\,E$ in (2.8) are the
same ones as in (2.3).

\medskip
{\bf Remark 2}. From Stirling's formula for the gamma-function it
follows that $R(\k_j) \ll \k_j^{-1/2}$. In view of (2.7) this
means that the series in (2.8) is absolutely convergent and
uniformly bounded in $s$ when $s = \s$ is real. Therefore, when $s
= \s \to 0+$, (2.8) gives a refinement of (2.6).

\medskip
{\bf Remark 3}. From (2.3) and (2.6) it transpires that
$\lambda(\s)$ is an error term when $0 < \s < 1$. For this reason
we considered the values $0 < |s| \le 1$ in (2.8), although one
could treat the case $|s| > 1$ as well.

\medskip
{\bf Remark 4}. From (2.8) and elementary properties of the
Laplace transform one can easily obtain the Laplace transform of
$$
E_2(T) := \int_0^T|\zt|^4\d t - TP_4(\log T),\quad P_4(x) =
\sum_{j=0}^4a_jx^j,\eqno(2.11)
$$
where $a_4 = 1/(2\pi^2)$ (for the evaluation of the remaining
coefficients $a_j$, see [7]).

\medskip
\section{Analytic continuation of the Mellin transform}

\medskip

Remarks on the general problem of analytic continuation of
the modified Mellin transform $\Z_k(s)$
were given in [11] and [17]. We start here by proving a general result, which
links the problem to the moments of $|\zt|$. This is a new result, which we state
as

\bigskip
THEOREM 2. {\it Let $k\in\NN$ be fixed. The bound
$$
\int_0^T|\zt|^{2k}\d t \;\ll_\e\; T^{c+\e}\eqno(3.1)
$$
holds for some constant $c$, if and only if $\Z_k(s)$ is regular for
$\R s > c$, and for any given $\e>0$}
$$
\Z_k(c+\e+it) \;\ll_\e\; 1.\eqno(3.2)
$$

\bigskip
{\bf Proof of Theorem 2}. The constant $c$ must satisfy $c\ge1$ in view of
the known lower bounds for moments of $|\zt|$ (see e.g., [5, Chapter 9]).
Suppose that (3.1) holds. Then we have
$$
\int_X^{2X}|\zx|^{2k}x^{-s}\d x \;\ll
\; X^{-\s}\int_0^{2X}|\zx|^{2k}\d x
\ll_\e X^{c-\s+\e/2},
$$
where $\s = \R s$. Therefore
$$
\int_1^\infty |\zx|^{2k}x^{-s}\d x \;=\;
\sum_{j=0}^\infty \,\int_{2^j}^{2^{j+1}}|\zx|^{2k}x^{-s}\d x
\;\ll\; \sum_{j=0}^\infty 2^{j(c-\s+\e/2)} \;\ll_\e\; 1
$$
if $\s = c+\e$, since $\sum_j2^{-\e j/2}$ converges.
This shows that $\Z_k(s)$ is regular for $\s > c$ and that
(3.2) holds.

Conversely, suppose that $\Z_k(s)$ is regular for $\s > c$ and that
(3.2) holds. Using the classical integral
$$
{\rm e}^{-x} = {1\over2\pi i}\int_{(d)}x^{-s}\G(s)\d s\qquad
(\R x > 0,\,d>0),
$$
we have
\begin{eqnarray}\nonumber
\int_1^\infty {\rm e}^{-x/T}|\zx|^{2k}\d x
&=&
\int_1^\infty {1\over2\pi i}\int_{(c+\e)}\G(s)\bigl({x\over T}\bigr)^{-s}\d s
|\zx|^{2k}\d x\\\nonumber
&= &{1\over2\pi i}\int_{(c+\e)}\G(s) T^s\Z_k(s)\d s \ll_\e T^{c+\e},
\end{eqnarray}
by absolute convergence and the fast decay of the gamma-function. This
yields
\begin{eqnarray}\nonumber
\int_0^T|\zx|^{2k}\d x &\le& O(1) + {\rm e}\int_1^T
{\rm e}^{-x/T}|\zx|^{2k}\d x\\ \nonumber
&\ll& 1 + \int_1^\infty {\rm e}^{-x/T}|\zx|^{2k}\d x
\ll_\e T^{c+\e},\nonumber
\end{eqnarray}
which proves (3.1).

\bigskip
{\bf Corollary 1}. The Lindel\"of hypothesis ($|\zt| \ll_\e |t|^\e$)
is equivalent to the statement that, for every $k\in\NN$, $\Z_k(s)$ is
regular for $\s>1$ and satisfies $\Z_k(1+\e+it) \ll_{k,\e} 1$.

\medskip
Indeed, the Lindel\"of hypothesis is equivalent (see e.g., [6, Section 1.9])
to (3.1) with $c=1$ for every $k\in\NN$. Therefore the assertion follows
from Theorem 1.

\bigskip
{\bf Corollary 2}. If we define
$$
\s_k := \inf\left\{\,d_k \;:\; \int_0^T|\zt|^{2k}\d t \ll_{k}
T^{d_k}\,\right\},\eqno(3.3)
$$
$$
\rho_k := \inf\left\{\,r_k \;:\; \Z_k(s) \;{\rm {is\, regular\, for}}\;
\R s > r_k\,\right\},\eqno(3.4)
$$
then
$$
\rho_k \;=\;\s_k,\quad \s_k \ge\;1.\eqno(3.5)
$$

\medskip
Note that from the  bounds on $I_k(T)$ in [5, Chapter 8] we obtain
$$
\s_1 = \s_2 = 1,\quad \s_k \le {k+2\over4}\quad(3\le k \le 6),\eqno(3.6)
$$
and upper bounds for $\s_k$ when $k>7$ may be obtained by using results
on the corresponding higher power moments of $|\zt|$ (op. cit.).
The Lindel\"of hypothesis may be reformulated as $\s_k = 1\;(\forall k \ge 1)$.

\smallskip
Thus at present we have two situations regarding analytic continuation
of $\Z_k(s)$:

\smallskip
a) For $k = 1,2,$ one can obtain analytic continuation of $\Z_k(s)$ to
the left of $\R s = 1$ (in fact to $\CC$). This will be discussed in
Section 5 and Section 6, respectively.

\smallskip

b) For $k>2$ only upper bounds for $\s_k$ (cf. (3.6)) are known. A
challenging problem is to improve these bounds, which would entail
progress on bounds of power moments of $|\zt|$, one of the central
topics in the theory of $\z(s)$.

\medskip
In what concerns power moments of $|\zt|$ one expects to have a formula
analogous to (2.11). Namely
for any fixed $k \in\NN$, we expect
$$
\int_0^T|\zt|^{2k}\d t = TP_{k^2}(\log T) + E_k(T)\eqno(3.7)
$$
to hold, where it is generally assumed that
$$
P_{k^2}(y) = \sum_{j=0}^{k^2}a_{j,k}y^j\eqno(3.8)
$$
is a polynomial in $y$ of degree $k^2$ (the integral in (3.7) is
$\gg_k T\log^{k^2}T$; see e.g., [5, Chapter 9]). The function
$E_k(T)$ is to be considered as the error term in (2.7), namely
one supposes that
$$
E_k(T) \= o(T)\qquad(T \to \infty).\eqno(3.9)
$$
So far (3.7)--(3.9) are known to hold only for $k = 1$ and $k = 2$
(see [6] and [28] for a comprehensive account). Therefore in view of the
existing knowledge on the higher moments of $|\zt|$,
embodied in (3.6), at present
the really important cases of (3.7) are $k = 1$ and $k = 2$. Plausible
heuristic arguments for the values of the coefficients $a_{j,k}$ were
recently given by Conrey et al. [3], by using methods from Random Matrix Theory
(see also Keating--Snaith [23]).

\medskip
In case (3.7)--(3.9) hold, this may be used to obtain the
analytic continuation of $\Z_k(s)$ to the region $\s\ge 1$ (at least).
Indeed, by using (3.7)-(3.9) we have
$$
\Z_k(s) = \int_1^\infty |\zx|^{2k}x^{-s}\d x = \int_1^\infty
x^{-s}\d\left(xP_{k^2}(\log x) + E_k(x)\right)
$$ $$
= \int_1^\infty (P_{k^2}(\log x) + P'_{k^2}(\log x))x^{-s}\d x
- E_k(1) + s\int_1^\infty E_k(x)x^{-s-1}\d x.\eqno(3.10)
$$
But for $\R s > 1$ change of variable $\log x = t$ gives
$$
\int_1^\infty (P_{k^2}(\log x) + P'_{k^2}(\log x))x^{-s}\d x
 = \int_1^\infty \left\{\sum_{j=0}^{k^2}a_{j,k}\log^jx +
\sum_{j=0}^{k^2-1}(j+1)a_{j+1,k}\log^jx\right\}x^{-s}\d x
$$
$$
= \int_0^\infty \left\{\sum_{j=0}^{k^2}a_{j,k}t^j +
\sum_{j=0}^{k^2-1}(j+1)a_{j+1,k}t^j\right\}{\rm e}^{-(s-1)t}\d t
= {a_{k^2,k}(k^2)!\over(s-1)^{k^2+1}} +
\sum_{j=0}^{k^2-1} {a_{j,k}j! + a_{j+1,k}(j+1)!\over
(s-1)^{j+1}}.
\eqno(3.11)
$$
Hence inserting (3.11) in (3.10) and using (3.9)
we obtain  the analytic continuation of $\Z_k(s)$
to the region $\s\ge1$. As we know (e.g., see [6] and [28]) that
$$
\int_1^T E_1^2(t)\d t \ll T^{3/2},\qquad
\int_1^T E_2^2(t)\d t \ll T^{2}\log^{22}T,\eqno(3.12)
$$
it follows on applying the Cauchy--Schwarz inequality to the last
integral in (3.10) that (3.9)-(3.11) actually provides
the analytic continuation of $\Z_1(s)$ to the region $\R s > 1/4$, and of
$\Z_2(s)$ to $\R s > 1/2$.

\medskip
\section{Recurrence formulas and identities}
\medskip
There is a possibility to obtain analytic continuation
of $\Z_k(s)$ by using a recurrent relation involving
$\Z_r(s)$ with $r < k$, which was mentioned in [11] and [17].
This result of ours is

\medskip
THEOREM 3. {\it For $k\ge 2,\,r = 1,\ldots, k-1$, $\R s$ and $c = c(k,r)$ sufficiently
large, we have}
$$
{\cal Z}_k(s) = {1\over2\pi i}\int_{(c)}{\cal Z}_{k-r}(w){\cal Z}_r(1+s-w)
\d w.\eqno(4.1)
$$

\medskip
{\bf Proof of Theorem 3}. For $\R(1-s)$ sufficiently large we have
$$
{\cal Z}_k(1-s) = \int_0^\infty \zeta^*(x)x^{s-1}\d x,
$$
where $\zeta^*(x) = |\zx|^{2k}$ if $x \ge 1$ and zero otherwise. If $F(s)$ is
the Mellin transform of $f(x)$ (see e.g., the Appendix of [5] for conditions
under which this holds) then one has the
Mellin inversion formula
$$
\hf(f(x+0) + f(x-0)) = {1\over2\pi i}\lim_{T\to\infty}\int_{\s-iT}^{\s+iT}
F(s)x^{-s}\d s.\eqno(4.2)
$$
The use of this relation gives
$$
|\zx|^{2k} = {1\over2\pi i}\int_{(c)}{\cal Z}_k(1-s)x^{-s}\d s,
\qquad(c \le c_0(k) < 0,\;x \ge 1).\eqno(4.3)
$$
Therefore, for $\,k,r \in \NN,\,k\ge2,\,1 \le r < k$, by using (4.3) we obtain
$$
\int_1^\infty |\zx|^{2k}x^{-s}\d x = \int_1^\infty |\zx|^{2r}
|\zx|^{2(k-r)}x^{-s}\d x
$$
$$
= \int_1^\infty |\zx|^{2r}\left({1\over2\pi i}\int_{(c)}{\cal Z}_{k-r}
(1-w)x^{-w}\d w\right)x^{-s}\d x
$$
$$
= {1\over2\pi i}\int_{(c)}{\cal Z}_{r}(w+s){\cal Z}_{k-r}(1-w)\d w
\quad(\s \ge \s_0(k) > 1-c).
$$
Changing $1-w$ to $w$ we obtain (4.1).

\medskip
In particular, by using (3.6), we obtain the identities
$$
\Z_3(s) \= {1\over2\pi i}\int_{(1+\e)}{\Z}_1(w)
\Z_2(1 - w + s)\d w\qquad(\s > {\txt{5\over4}}),
$$
$$
{\Z}_4(s) \= {1\over2\pi i}\int_{({5\over4}+\e)}{\Z}_2(w)
{\Z}_2(1 - w + s)\d w\qquad(\s > {\txt{3\over2}}).
$$

\medskip

The following result provides an integral representation for $\Z_k^2(s)$.
This is

\bigskip
THEOREM 4. {\it In the region of absolute convergence we have}
$$
\Z_k^2(s) = 2\int_1^\infty x^{-s}\left(\int_{\sqrt{x}}^x|\z(\hf + iu)|^{2k}
\Bigl|\z\Bigl(\hf +i\,{x\over u}\Bigr)\Bigr|^{2k}{\d u\over u}\right)\d x.
\eqno(4.4)
$$

\bigskip
{\bf Proof of Theorem 4}. Clearly the region of validity
of (4.4) depends on $k$, and by using e.g. (3.6) one can
provide explicit $\bar{\s}_k$ such that (4.4) holds
for $\s > \bar{\s}_k$. To prove the assertion,
we set $f(x) = |\zx|^{2k}$ and make the change
of variables $xy = X,\,x/y = Y$, so that the absolute value of the
Jacobian of the transformation is equal to $1/(2Y)$. Therefore
$$
\Z_k^2(s) = \int_1^\infty\int_1^\infty (xy)^{-s}f(x)f(y)\d x\d y
= {1\over2}\int_1^\infty X^{-s}\int_{1/X}^X{1\over Y}
f(\sqrt{XY}\,)f(\sqrt{X/Y}\,)\d Y\d X.
$$
But as we have ($y = 1/u$)
$$
\int_{1/x}^x
f(\sqrt{xy}\,)f(\sqrt{x/y}\,){\d y\over y} = \int_1^xf(\sqrt{x/u})
f(\sqrt{xu}\,){\d u\over u},
$$
we obtain that, in the region of absolute convergence, the identity
$$
\Z_k^2(s) = \hf\int_1^\infty x^{-s}\left(\int_1^x f(\sqrt{xy}\,)
f(\sqrt{x/y}\,){\d y\over y}\right)\d x
$$
is valid. The inner integral here becomes, after the change of
variable $\sqrt{xy} = u$,
$$
2\int_{\sqrt{x}}^x f(u)f\bigl({x\over u}\bigr){\d u\over u},
$$
and (4.4) follows. The argument also shows that, for  $0 < a < b$
and any integrable function $f$ on $[a,\,b]$,
$$
\left(\int_a^b f(x)x^{-s}\d x\right)^2
= 2\int_{a^2}^{b^2}x^{-s}\left\{\int_{\sqrt{x}}^{\min(x/a,b)}f(u)
f\bigl({x\over u}\bigr){\d u\over u}\right\}\d x.
$$

\medskip
\section{The modified Mellin transform, $k=1$}

\medskip

We begin our discussion of
the function $\Z_1(s)$ by obtaining its analytic continuation over $\CC$.
The relevant result is contained in

\bigskip
{THEOREM 4}. {\it The function $\Z_1(s)$
continues meromorphically to $\CC$, having only a double pole
at $s=1$, and at most simple poles at $s = -1,-3,\ldots\,$. The principal
part of its Laurent expansion at $s=1$ is given by
$$
{1\over(s-1)^2} + {2\gamma - \log(2\pi)\over s - 1},\eqno(5.1)
$$
where $\gamma $ is Euler's constant.}

\bigskip
{\bf Proof of Theorem 4}. It was shown in [17] that $\Z_1(s)$
continues analytically to a function that is regular for $\s > -3/4$. In [22]
M. Jutila proved
that $\Z_1(s)$ continues meromorphically to $\CC$, having only a double pole
at $s=1$ and at most double poles for $s = -1,-2,\ldots\,$. The present form
of Theorem 4 was obtained by the author in [13], and a different proof
is to be found in the dissertation of M. Lukkarinen [25].
As our proof of Theorem 4 is fairly simple, it will be given now. Let
$$
{\bar L}_1(s) := \int_1^\infty|\z(\hf+iy)|^{2}{\rm e}^{-ys}
\d y\quad( \R s > 0). \eqno(5.2)
$$
Then we have by absolute convergence, taking $\s = \R s$ sufficiently large
and making the change of variable $xy = t$,
\begin{eqnarray*}
\int_0^\infty {\bar L}_1(x) x^{s-1}\d x &=& \int_0^\infty
\left(\int_1^\infty|\z(\hf + iy)|^{2}{\rm e}^{-yx}\d y\right)
x^{s-1}\d x\\
&=& \int_1^\infty|\z(\hf + iy)|^{2}\left(\int_0^\infty
x^{s-1}{\rm e}^{-xy}\d x\right)\d y
\end{eqnarray*}
$$
= \int_1^\infty|\z(\hf + iy)|^{2}y^{-s}\d y \int_0^\infty
{\rm e}^{-t}t^{s-1}\d t = {\cal Z}_1(s)\G(s).\eqno(5.3)
$$
Further we have
\begin{eqnarray*}
\int_0^\infty {\bar L}_1(x) x^{s-1}\d x &=& \int_0^1{\bar L}_1(x) x^{s-1}\d x
+ \int_1^\infty{\bar L}_1(x) x^{s-1}\d x\\
&=& \int_1^\infty{\bar L}_1(1/x) x^{-1-s}\d x + A(s)\quad(\s>1),
\end{eqnarray*}
say, where $A(s)$ is an entire function. Since (see (1.1))
$$
{\bar L}_1(1/x) =  L_1(1/x) - \int_0^1|\z(\hf+iy)|^2{\rm e}^{-y/x}
\d y\qquad(x\ge1),
$$
it follows from (5.3) by analytic continuation that, for $\s>1$,
$$
{\cal Z}_1(s)\G(s) = \int_1^\infty L_1(1/x) x^{-1-s}\d x
- \int_1^\infty\Bigl(\int_0^1|\z(\hf+iy)|^2{\rm e}^{-y/x}
\d y\Bigr)x^{-1-s}\d x + A(s)
$$ $$
= I_1(s) - I_2(s) + A(s),\eqno(5.4)
$$
say. Clearly, for any integer $M\ge1$, we have
$$
I_2(s)= \int_1^\infty\int_0^1|\z(\hf+iy)|^2
\left(\sum_{m=0}^M{(-1)^m\over m!}\left({y\over x}\right)^m +
O_M(x^{-M-1})\right)
\d y \,x^{-1-s}\d x$$
$$
= \sum_{m=0}^M{(-1)^m\over m!}h_m\cdot{1\over m+s} + H_M(s),\eqno(5.5)
$$
say, where $H_M(s)$ is a regular function of $s$ for $\s > -M-1$, and
$h_m$ is a constant. Note that, for $\s = 1/T\,(T\to\infty)$
and any $N\ge0$, (2.1) gives
$$
L_1\left({1\over T}\right) =
\left(\log\left({T\over2\pi}\right) +
\gamma\right)\sum_{n=0}^N a_nT^{1-2n}
+ \sum_{n=0}^N b_nT^{-2n} + O_N(T^{-1-2N}\log T)
$$
with suitable $a_n, b_n\,(a_0 = 1)$. Inserting this formula in $I_1(s)$ in (5.4)
we have
$$
I_1(s) = \int_1^\infty (\log {x\over2\pi}+\gamma)
\sum_{n=0}^Na_nx^{-2n-s}\d x + \int_1^\infty \sum_{n=0}^Nb_nx^{-1-2n-s}\d x
+ K_N(s)$$ $$
= \sum_{n=0}^Na_n\left({1\over(2n+s-1)^2} +
{\gamma-\log2\pi\over 2n+s-1}\right)
+  K_N(s)\quad(\s > 1),\eqno(5.6)
$$
say, where $K_N(s)$ is regular for $\s > -2N$. Taking $M = 2N$
it follows from (5.4)--(5.6) that
$$
\Z_1(s)\G(s) = \sum_{n=0}^Na_n\left({1\over(2n+s-1)^2} +
{\gamma-\log2\pi\over 2n+s-1}\right)
+ \sum_{m=0}^{2N}{(-1)^m\over m!}h_m\cdot{1\over m+s} + R_N(s),
\eqno(5.7)
$$
say, where $R_N(s)$ is a regular function of $s$ for $\s > -2N$. This holds
initially for $\s >1$, but by analytic continuation it holds for $\s > -2N$.
Since $N$ is arbitrary
and $\G(s)$ has no zeros, it follows that (5.7)
provides meromorphic continuation
of $\Z_1(s)$ to $\CC$. Taking into
account that $\G(s)$ has simple poles at $s = -m\;(m=0,1,2,\ldots\,)$
we obtain then the analytic
continuation of $\Z_1(s)$ to $\CC$, showing that besides $s=1$
the only poles of $\Z_1(s)$ can be
simple poles at $s = 1-2n$ for $n \in\NN$, as asserted by Theorem 4.
With more care the residues at these poles could be explicitly evaluated.
Finally using (5.7) and
$$
{1\over \G(s)} = 1 + \gamma(s-1) + \sum_{n=2}^\infty d_n(s-1)^n
$$
we obtain that the principal part of the Laurent expansion at $s=1$ is
given by (5.1).

\bigskip
Another major problem is to determine the order of growth of $\Z_1(\s+it)$.
Concerning pointwise bounds  of $\Z_1(s)$, we have (see M. Jutila [22])
$$
\Z_1(\s+it) \ll_\e t^{{5\over6}-\s+\e}\qquad(\hf \le \s \le 1,\;t\ge t_0).
\eqno(5.9)
$$
We also have the mean square bounds (see [17] for proof)
$$
\int_1^T|{\cal Z}_1 (\s+it)|^2\d t \ll _\e T^{3-4\s+\e}\qquad(
0 \le \s \le{\txt{1\over2}}),\eqno(5.10)
$$
and
$$
\int_1^T|{\cal Z}_1 (\s+it)|^2\d t \ll _\e
T^{2-2\s+\e}\qquad({\txt{1\over2}} \le \s \le 1).
\eqno(5.11)
$$
The bound in (5.11) is essentially best possible
since, for any given $\e > 0$,
$$
\int_1^T|\Z_k(\s + it)|^2\d t \;\gg_\e\; T^{2-2\s-\e}\qquad(k = 1,2;\;
\hf < \s < 1).\eqno(5.12)
$$
This assertion follows from
$$
\int_T^{2T}|\zt|^{2k}\d t \ll_\e
T^{2\s-1}\int_0^{T^{1+\e}}|\Z_k(\s+it)|^2\d t \quad(k = 1,2;\,\hf < \s < 1)
\eqno(5.13)
$$
and lower bounds for the integral on the left-hand side
(see [5, Chapter 9]).
The proof of (5.13) when $k = 2$  appeared in [11], and the
proof of the bound when $k=1$ is on similar lines.
It is plausible to conjecture that (5.9) holds with the exponent 1/2 instead
of 5/6 on the right-hand side.

\medskip
\section{The modified Mellin transform, $k=2$}

\medskip

The function ${\cal Z}_2(s)$ has quite a different analytic behaviour
from the function $\Z_1(s)$. It
was introduced and studied by Y. Motohashi [26]
(see also his monograph [28]). He  has shown   that ${\cal Z}_2(s)$ has meromorphic
continuation over $\CC$. In the half-plane $\R s > 0$ it has
the following singularities: the pole $s = 1$ of order five,  simple
poles at $s = {1\over2} \pm i\k_j\,(\k_j =\sqrt{\lambda_j -
{1\over4}})$ and poles at $s =
\hf\rho$,  where $\rho$ denotes complex zeros of
$\zeta(s)$. The residue of ${\cal Z}_2(s)$ at
$s = {1\over2} + i\k_h$ equals (see Section 2 for definitions concerning
spectral theory)
$$ R(\k_h) := \sqrt{{\pi\over2}}{\Bigl(2^{-i\k_h}{\G({1\over4} -
\hf{i}\k_h)\over\G({1\over4} +
\hf{i}\k_h)}\Bigr)}^3\G(2i\k_h)\cosh(\pi\k_h)
\sum_{\k_j=\k_h}\alpha_j H_j^3(\txt{1\over2}),
$$
and the residue at $s = {1\over2} - i\k_h$ equals
$\overline{R(\k_h)}$.
The principal part of ${\cal Z}_2(s)$ has the form (this
may be compared with (5.1))
$$
\sum_{j=1}^5{A_j\over(s-1)^j},\eqno(6.1)
$$
where $A_5 = 12/\pi^2$, and the remaining $A_j$'s can be
evaluated explicitly by following the analysis in [26].
The function ${\cal Z}_2(s)$ was used to furnish several strong
results on $E_2(T)$ (see (2.7)), the error term in the asymptotic
formula for the fourth moment of $|\zt|$. Y. Motohashi [26]
used it to show that $E_2(T) = \Omega_\pm(T^{1/2})$,
which sharpens the earlier result of Ivi\'c-Motohashi (see [18])
that $E_2(T) = \Omega(T^{1/2})$.
The same authors (see e.g., [19], [20] and [28]) have proved that
$$
E_2(T)  \ll T^{2/3}\log^{C_1}T\;(C_1>0),
\quad\int_0^TE_2(t)\d t \ll T^{3/2},\eqno(6.2)
$$
as well as the second bound in (2.12).  In [8] and [9] the
author has applied the theory of $\Z_2(s)$ to obtain the following
quantitative omega-results:
There exist constants $A, B > 0$ such that for $T \ge T_0 > 0$ every
interval $\,[T,\,AT]\,$ contains points $t_1, t_2, t_3, t_4$ such that
$$
E_2(t_1) > Bt_1^{1/2},\; E_2(t_2) < -Bt_2^{1/2},\;
\int_0^{t_3}E_2(t)\d t > Bt_3^{3/2},\;
\int_0^{t_4}E_2(t)\d t < -Bt_4^{3/2}.
$$
Moreover, we have (see [9])
$$
\int_0^TE_2^2(t)\d t \gg T^2,\eqno(6.3)
$$
which complements the upper bound in (6.2).

\medskip
As for the estimation of $\Z_2(s)$, we have (see the author's work [12])
$$
\Z_2(\s + it) \;\ll_\e\; t^{{4\over3}(1-\s)+\e}\qquad(\hf < \s \le 1;\,
t \ge t_0 > 0),\eqno(6.4)
$$
and (similarly to (5.9)) I have conjectured in [12] that the exponent on
the right-hand side of (6.4) can be replaced by $1/2 - \s$. This is very
strong, as it implies (6.12) with $\theta = 1$ and $E_2(T) \ll_\e T^{1/2+\e}$,
and both bounds (up to ``$\e$") are best possible.
It was proved in [17] that
$$
\int_0^T|{\cal Z}_2(\s + it)|^2\d t \ll_\e
T^\e\left(T + T^{2-2\s\over1-c}\right) \qquad(\hf < \s \le 1),\eqno(6.5)
$$
and we also have unconditionally
$$
\int_0^T|\Z_2(\s + it)|^2\d t \;\ll\;T^{10-8\s\over3}\log^CT
\qquad(\hf < \s \le 1,\,C > 0).\eqno(6.6)
$$
The constant $c$ appearing in (6.5) is defined by $
E_2(T) \ll_\e T^{c+\e},$ so that by (6.2) and (6.3) we have
$\hf \le c \le {2\over3}$. In (6.4)--(6.6) $\s$ is assumed
to be fixed, as $s = \s+it$ has to stay away from the line $\R s = \hf$
where $\Z_2(s)$ has poles. Lastly, the author [14] proved that,
for ${5\over6} \le \s \le {5\over4}$ we have,
$$
\int_1^{T}|\Z_2(\s+it)|^2\d t \ll_\e T^{{15-12\s\over5}+\e}.\eqno(6.7)
$$
The lower limit of integration in (6.7) is unity, because of the pole
$s=1$ of $\Z_2(s)$. By taking $c = 2/3$ in (6.5) and using the convexity
of mean values (see [5, Lemma 8.3]) it follows that
$$
\int_1^T|\Z_2(\s+it)|^2\d t \;\ll_\e\; T^{{7-6\s\over2}+\e}\qquad(\hf < \s
\le \txt{5\over6}).\eqno(6.8)
$$
Note that (6.7) and (6.8) combined  provide the sharpest known bounds in the
whole range $\hf < \s \le {5\over6}$.

\medskip
Both pointwise and mean square estimates for $\Z_2(s)$ may be used
to estimate $E_2(T)$ and the eighth moment of $|\zt|$.
This connection is furnished by the following
result, proved by the author in [12].

\bigskip
THEOREM 6. {\it Suppose that, for some $\rho \ge 0$ and $r \ge 0$,}
$$
\Z_2(\s + it) \;\ll_\e\; |t|^{\rho+\e},\quad
\int_1^T|\Z_2(\s + it)|^2\d t \;\ll_\e\; T^{1+2r+\e}\quad(\hf < \s \le 1),
\eqno(6.9)
$$
{\it where $\s$ is fixed and $|t|\ge t_0 > 0$. Then we have}
$$
E_2(T) \;\ll_\e\; T^{{2\rho+1\over2\rho+2}+\e}, \quad
E_2(T) \;\ll_\e\; T^{{2r+1\over2r+2}+\e}\eqno(6.10)
$$
{\it and}
$$
\int_0^T|\zt|^8\d t \;\ll_\e\; T^{{4r+1\over2r+1}+\e}.\eqno(6.11)
$$
\medskip\no
Note that the conditions $\rho\ge0,\,r\ge 0$ must hold in view of (5.12).
Also note that from (6.6) with $\s = \hf + \e$ one can take in
(6.9) $r = \hf$, hence (6.10) gives
$E_2(T) \;\ll_\e\; T^{{2\over3}+\e},$
which is essentially the strongest known bound (see (6.2)). Thus
any improvement of the existing mean square bound for $\Z_2(s)$ at
$\s = \hf + \e$ would result in the bound for $E_2(T)$
 with the exponent strictly less than
2/3, which would be important.  Of course, if the first bound in (6.9)
holds with some $\rho$, then trivially the second bound will hold with
$r = \rho$, i.e. $r \le \rho$ has to hold. Observe that the
known value $r = \hf$ and (6.11) yield
$$
\int_0^T|\zt|^8\d t \;\ll_\e\; T^{\theta+\e}\eqno(6.12)
$$
with $\theta = 3/2$, which is, up to``$\e$", currently
the best known upper bound
(see [5, Chapter 8]) for the eighth moment, and any value $r < \hf$
in (6.9) would reduce the exponent $\theta = 3/2$ in (6.12). The connections
between upper bounds for the integral in (6.12) and mean
square estimates involving $\Z_2(s)$ and related functions are
also given as

\bigskip
THEOREM 7. {\it The eighth moment bound, namely} (6.12) {\it with
$\theta = 1$, is equivalent to the mean square bound
$$
\int_1^T|\Z_2(1+it)|^2\d t \;\ll_\e\; T^{\e},\eqno(6.13)
$$
and to
$$
\int_T^{2T}I^2(t,G)\d t \;\ll_\e\; T^{1+\e}\qquad
(T^\e \le G = G(T) \le T),\eqno(6.14)
$$
where}
$$
\quad I(T,G) := {1\over\sqrt{\pi}G}\int_{-\infty}^\infty
|\z(\hf + iT + iu)|^4{\rm e}^{-(u/G)^2}\d u.\eqno(6.15)
$$

\bigskip
{\bf Proof of Theorem 7}. We suppose first that (6.13) holds. Then
(6.12) with $\theta = 1$ follows from (5.13) with $k=2$.
Conversely, if (6.12) holds with $\theta = 1$,
note that we have, by [11, Lemma 4],
$$
\int_T^{2T}\left|\int_X^{2X}|\zx|^4x^{-s}\d x\right|^2\d t
\ll \int_X^{2X}|\zx|^8x^{1-2\s}\d x
\ll_\e X^{2-2\s+\e}\eqno(6.16)
$$
for $s = \s +it,\;\hf < \s \le 1$.
Similarly, using the Cauchy-Schwarz inequality for integrals
and the second bound in (3.12), it follows that
$$
\int_T^{2T}\left|\int_X^{2X}E_2(x)x^{-s}\d x\right|^2\d t
\ll_\e X^{1-2\s+\e}\qquad(s = \s + it,\;\s > \hf).\eqno(6.17)
$$
Combining (6.16) and (6.17) we obtain then,
similarly to the proof of (5.7) in [17],
$$
\int_1^T|\Z_2(\s+it)|^2\d t \;\ll_\e\; T^{4-4\s+\e}\qquad(\hf < \s\le 1),
$$
and for $\s=1$ we have (6.13). From (6.7) it follows that the integral in (6.13)
is unconditionally bounded by $T^{3/5+\e}$, and any improvement of the exponent
3/5 would also result in the improvement of the exponent
$\theta = 3/2$ in (6.12).

\medskip
Suppose again that (6.12) holds with $\theta = 1$. Then the left-hand
side of (6.14) is, for $T^\e \le G = G(T) \le T$,
\begin{eqnarray*}
 &&\int_T^{2T}\left({1\over\sqrt{\pi}G}\int_{-\infty}^\infty
|\z(\hf + it + iu)|^4{\rm e}^{-(u/G)^2}\d u\right)^2\d t\\
&\ll& 1 + G^{-2}\int_T^{2T}\left(\int_{-G\log T}^{G\log T}
|\z(\hf + it + iu)|^4{\rm e}^{-(u/G)^2}\d u\right)^2\d t\\
&\ll& 1 + G^{-1}\int_{-G\log T}^{G\log T}
\left(\int_T^{2T}|\z(\hf + it + iu)|^8\d t\right)\d u\\
&\ll_\e& G^{-1}\int_{-G\log T}^{G\log T} T^{1+\e}\d u\\
&\ll_\e&  T^{1+\e},
\end{eqnarray*}
as asserted. We remark that (6.14) is trivial when $T^{2/3} \le G \le T$
(see e.g., [6, Chapter 5]). Finally, if (6.14) holds, then we use
[15, Theorem 4], which in particular says that, for fixed $m\in\NN$,
$$
\int_T^{2T}I^m(t,G)\d t \ll_\e T^{1+\e}
\quad(T^{\a_{m}+\e}\le G = G(T) \le T,\,
 0 \le \a_{m} < 1)\eqno(6.18)
$$
implies that
$$
\int_0^T|\zt|^{4m}\d t \ll_\e T^{1+(m-1)\a_{m}+\e}.
$$
Using this result with $m = 2, \a_{2} = 0$, we obtain at once
(6.12) with $\theta = 1$. This completes the proof of Theorem 7.
So far (6.18) is known to hold unconditionally with $\a_2 = \hf$
(see [14]), which yields another proof of (6.12) with $\theta = 3/2$.

\bigskip
The significance of (6.14) is that for $I(T,G)$ in (6.15) an
explicit formula of Y. Motohashi (see [6], [26] or [28]) exists. It
involves quantities from spectral theory, and thus the eighth
moment problem is directly connected to this theory via Theorem 7.
Namely
$$
I(T,G)= {\pi\over\sqrt{2T}}\sum_{j=1}^\infty \a_j
H_j^3(\hf) \k_j^{-{1\over2}}
\sin\left(\k_j\log{\k_j\over 4{\rm e}T}\right)
\times\exp\Biggl(-{1\over4}
\Bigl({G\k_j\over T}\Bigr)^2\Biggr) + O(\log^{3D+9}T),
$$
if $\;T^{1/2}\log^{-D}T \le G \le T/\log T$ for an arbitrary, fixed
constant $D>0$. We have the following new result, proved by the author in [16]:

\smallskip
THEOREM 8. {\it Let, for $m\in\NN$ and
$1 \ll K < K'\le 2K\ll T, T \le t \le 2T$,
$$
S_m(K;K',t) =
\sum_{K<\k_j\le K'<2K} \a_j H_j^m(\hf)\cos\Bigl(\k_j\log\bigl(
{4{\rm e}t\over \k_j}\bigr)\Bigr).
$$
Then, for $m = 1,2,3$,}
$$
\int_T^{2T}{\bigl(S_m(K;K',t)\bigr)}^2\d t \,\ll_\e\, T^{1+\e}K^3.
$$

\medskip
{\bf Corollary.} {\it We have}
$$
\int\limits_0^T E_2^2(t)\d t \;\ll_\e\; T^{2+\e},\quad
\int\limits_0^T|\zt|^{12}\d t \;\ll_\e\; T^{2+\e}.\eqno(6.19)
$$

\smallskip\no
The first bound (see [9], [19]), up to ``$\e$'', is best possible, the second
is a consequence of the first and is a well-known result of D.R. Heath-Brown [4].
We shall just give a sketch of the first bound in (6.19). From [6, Lemma 5.1]
we have
$$
E_2(2T) - E_2(T) \le S(2T+\D\log T,\D) - S(T-\D\log T,\D)
+  O(\D\log^5T) + O(T^{1/2}\log^CT)
$$
with $T^{1/2} \le\D\le T^{1-\e}$ and
$$
S(T,\D) := \pi\sqrt{\hf T}\sum_{j=1}^\infty \a_j\H \k_j^{-{3\over2}}
\cos\left(\k_j\log{\k_j\over 4{\rm e}T}\right)\exp\Bigl(-{\txt{1\over4}}
\bigl({\D\k_j\over T}\bigr)^2\Bigr).
$$
We can truncate the above series at $T\D^{-1}\log T$ with a negligible error
and remove, by partial summation, the monotonic coefficients
$\k_j^{-3/2}$ and $\exp\Bigl(-{\txt{1\over4}}
\bigl({\D\k_j\over T}\bigr)^2\Bigr)$.Then we
obtain  the sum $S_m(K;K',t)$ with $m=3$ and $t$ replaced by
$2t+\D\log T$ or $t-\D\log T$, which does not cause any trouble. We have
$$
\int_T^{2T}{(E_2(2t)-E_2(t))}^2\d t
\le \int_{T/2}^{5T/2}\f(t){(E_2(2t)-E_2(t))}^2\d t,\eqno(6.20)
$$
where $\f(t)\,(\ge0)$ is a smooth function suported in $[T/2,\,5T/2]$ which
equals unity in $[T,\,2T]$. Hence the
integral on the left-hand side of (6.20) is
essentially majorized by $\ll_\e T^\e$ integrals of the type
$$
T\int_{T/2}^{5T/2}\f(t){\bigl( K^{-3/2}S_m(K;K',t)\bigr)}^2\d t
\ll_\e T^{2+\e},
$$
with $m=3$ and (6.19) follows from Theorem 8 on
replacing $t$ by $t2^{-j}$ in the integrand in
(6.20), and summing up the corresponding bounds over $j = 1,2,\ldots\,$.



\small\bigskip
Aleksandar Ivi\'c

Katedra Matematike RGF-a

Universitet u Beogradu, \DJ u\v sina 7

11000 Beograd, Serbia

\tt ivic@rgf.bg.ac.yu, \enskip aivic@matf.bg.ac.yu
\end{document}